\documentclass[12pt]{amsart}
\usepackage[dvips]{epsfig}
\usepackage{graphics}
\usepackage{latexsym}
\usepackage{verbatim}
\usepackage{amsmath}
\usepackage{amsthm}
\usepackage{amssymb}
\usepackage{hyperref}
\newcommand{\ds}{\displaystyle}
\newcommand{\vip}{\vskip0.2cm}
\newcommand{\ala}{\nonumber \\}
\newcommand{\indiq}{1\!\! 1}
\newcommand{\e}{{\varepsilon}}

\newcommand{\ch}{{{\mathcal H}}}
\newcommand{\cp}{{{\mathcal P}}}

\newcommand{\rr}{{\mathbb{R}}}
\newcommand{\R}{{\mathbb{R}}}
\newcommand{\Sp}{{\mathbb{S}}}

\newcommand{\rd}{{\mathbb{R}^d}}
\newcommand{\rdd}{{\rd\times\rd}}
\newcommand{\lip}{{\mbox{Lip}(\rd)}}
\newcommand{\lipu}{{\mbox{Lip}_1(\rd)}}

\newcommand{\pun}{{{\mathcal P}_1(\rd)}}
\newcommand{\li}{L^\infty([0,T],\pun)}

\newcommand{\lilp}{L^\infty([0,T],\pun) \cap L^1([0,T], L^p(\rd))}

\newcommand{\intot}{\ds\int_0^t}
\newcommand{\intrd}{\ds\int_{\rd} \!\!\!}
\newcommand{\intrdd}{\ds\int_{\rd\times\rd}\!\!\!\!\!\!\!\!\!}
\newcommand{\intzp}{\ds\int_0^\pi}

\newcommand{\cost}{\cos \theta}
\newcommand{\sint}{\sin \theta}
\newcommand{\tf}{{\tilde f}}
\newcommand{\tg}{{\tilde g}}
\newcommand{\tv}{{\tilde v}}
\newcommand{\tz}{{\tilde z}}

\newcommand{\vs}{{v_*}}
\newcommand{\tvs}{{\tv_*}}

\renewcommand{\theequation}{\thesection.\arabic{equation}}

\newtheorem{theo}{\indent Theorem}[section]

\newtheorem{lem}[theo]{\indent Lemma}
\newtheorem{defin}[theo]{\indent Definition}
\newtheorem{cor}[theo]{\indent Corollary}

\begin{document}

\title[Well-posedness of some singular Boltzmann equations]
{On the well-posedness of the spatially homogeneous 
Boltzmann equation with a moderate angular singularity}

\author{Nicolas Fournier$^1$, Cl\'ement Mouhot$^2$}

\footnotetext[1]{Centre de Math\'ematiques,
Facult\'e de Sciences et Technologies,
Universit\'e Paris~12, 61, avenue du G\'en\'eral de Gaulle, 94010 Cr\'eteil 
Cedex, France, e-mail: {\tt nicolas.fournier@univ-paris12.fr}}

\footnotetext[2]{Ceremade, Universit\'e Paris IX-Dauphine,
Place du Mar\'echal de Lattre de Tassigny, 75775 Paris, France, 
e-mail: \texttt{clement.mouhot@ceremade.dauphine.fr}}

\def\shortauthorname{Nicolas Fournier, Cl\'ement Mouhot}

\def\abstractname{Abstract}

\begin{abstract}
We prove an inequality on the Kantorovich-Rubinstein distance --which can 
be seen as a particular case of a Wasserstein metric--  
between two solutions of the spatially homogeneous Boltzmann
equation without angular cutoff, but with a moderate angular singularity.
Our method is in the spirit of~\cite{f}. 
We deduce some well-posedness and stability results in the physically relevant 
cases of hard and moderately soft potentials. 

\noindent In the case of hard potentials, we relax the regularity 
assumption of \cite{dm}, but we need 
stronger assumptions on the tail of the distribution (namely 
some exponential decay). We thus obtain the first uniqueness result
for measure initial data.

\noindent In the case of moderately soft 
potentials, we prove existence and uniqueness 
assuming only that the initial
datum has finite energy and entropy (for very moderately soft
potentials), plus sometimes an additionnal moment condition. 
We thus improve significantly on all previous results, where weighted 
Sobolev spaces were involved.
\end{abstract}

\maketitle

\smallskip


\tableofcontents


\section{Introduction} 
\setcounter{equation}{0}

\subsection{The Boltzmann equation}  
We consider a spatially homogeneous gas in dimension $d \ge 2$ 
modeled by the Boltzmann equation. 
Therefore the time-dependent density $f=f_t(v)$ of 
particles with velocity $v\in \rd$ solves
\begin{eqnarray} \label{be}
\partial_t f_t(v) = \intrd dv_* \int_{\Sp^{d-1}} d\sigma B(|v-v_*|,\theta)
\big[f_t(v')f_t(v'_*) -f_t(v)f_t(v_*)\big],
\end{eqnarray}
where 
$$
v'=\frac{v+v_*}{2} + \frac{|v-v_*|}{2}\sigma, \quad 
v'_*=\frac{v+v_*}{2} -\frac{|v-v_*|}{2}\sigma
$$ 
and $\theta$ is the so-called {\em deviation angle} defined by 
$\cos \theta = \frac{(v-v_*)}{|v-v_*|} \cdot \sigma$.

\vip

The {\em collision kernel} $B=B(|v-v_*|,\theta)=B(|v'-v'_*|,\theta)$ 
is given by physics 
and is related to the microscopic interaction between particles. 
In dimension $d=3$ 
it is related to the probabilistic {\em cross-section}  $\hat B$ 
of the distribution of possible 
outgoing velocities $v'$ and $v'_*$ arising from a collision with 
two particles with velocities 
$v$ and $v_*$, by the formula $B = |v-v_*| \, \hat B$. 
We refer to the review papers of 
Desvillettes \cite{desvillettes} and Villani \cite{villani} for more details. 

\vip

Conservation of mass, momentum and kinetic energy
hold at least formally for solutions to (\ref{be}), that is for all $t\geq 0$,
$$
\intrd  f_t(v) \, \varphi(v) \, dv 
= \intrd f_0(v) \, \varphi(v) \, dv, \qquad \varphi = 1, v, |v|^2 
$$
and we classically may assume without loss of generality that 
$\int_{\rd} f_0(v) \, dv=1$.

\subsection{Assumptions on the collision kernel} 
We shall assume that the collision kernel takes the form 
\renewcommand\theequation{{\bf A1}}
\begin{equation}
B(|v-v_*|,\theta)\, \sin^{d-2} \theta =\Phi(|v-v_*|) \, \beta(d \theta)
\end{equation}
\renewcommand\theequation{\thesection.\arabic{equation}}
for some function $\Phi: \R_+ \mapsto \R_+$ 
and some nonnegative measure $\beta$ on $(0,\pi]$.

\vip

In the case of an interaction potential $V(s)=1/r^s$ in dimension $d=3$, 
with $s\in (2,\infty)$, one has 
\begin{equation} \label{betapuissance}
\Phi(z) = \mbox{cst} \, z^{\gamma}, 
\quad \beta(\theta) \sim \mbox{cst} \, \theta^{-1-\nu}, \; \; \mbox{ with } 
\gamma=\frac{s-5}{s-1}, \quad \nu=\frac{2}{s-1}.
\end{equation}
On classically names {\em hard potentials} the case when 
$\gamma\in (0,1)$ ({\em i.e.}, $s>5$ in dimension $d=3$), 
{\em Maxwellian molecules} the case when $\gamma=0$ ({\em i.e.}, $s=5$ 
in dimension $d=3$), and {\em soft potentials} the case 
when $\gamma\in (-d,0)$ ({\em i.e.}, $s\in (2,5)$ in dimension $d=3$).

\vip

Let us emphasize that $\int_{0+}  \beta(d \theta) = + \infty$,
which expresses the affluence of {\em grazing collisions}, but
in any case, 
\begin{equation}\label{gccs}
\int_0^\pi \theta^2 \, \beta(d \theta) < + \infty.
\end{equation}

\vip

In this paper we shall deal with a moderate angular singularity, 
that is we shall 
assume that the collision kernel satisfies 
\renewcommand\theequation{{\bf A2}}
\begin{equation}
\kappa_1=\int_0^\pi \theta \, \beta(d\theta) < + \infty,
\end{equation}
which corresponds to $s\in (3,\infty)$ in~\eqref{betapuissance}).

\vip

We will also assume that $\Phi$ behaves as a power function, 
namely that for some $\gamma \in (-d,1]$, there exists some
constant $C$ such that for all $z,\tz \in \R_+$,

\renewcommand\theequation{{\bf A3}($\gamma$)}
\begin{equation}
\Phi(z) \leq C \, z^\gamma; \quad |\Phi(z)- \Phi(\tz)|\leq C \, |z^\gamma -\tz^\gamma|.
\end{equation}

Sometimes, we will need a lowerbound: 
there exists $c>0$ such that for all $z \in \R_+$,

\renewcommand\theequation{{\bf A4}($\gamma$)}
\begin{equation}
\Phi(z) \geq c \, z^\gamma.
\end{equation}

In the case of hard potentials, we will also sometimes
use an additionnal technical assumption
in order to obtain the propagation of some exponential moments:

\renewcommand\theequation{{\bf A5}}
\begin{equation}
\beta(d\theta)=b(\cos \theta) \, d\theta, \mbox{ where }
b \mbox{ is nondecreasing, convex and } C^1 \ \mbox{ on } [-1,1).
\end{equation}

In the case of moderately soft potentials, we will sometimes use

\renewcommand\theequation{{\bf A6}($\nu$) }
\begin{equation} 
\beta(d\theta)=\beta(\theta) \, d\theta \hbox{ with }
\beta(\theta) \, \sim_{\theta \to 0} \mbox{cst} \, \theta^{-1-\nu}
\end{equation}
for some positive constant.
\renewcommand\theequation{\thesection.\arabic{equation}}

In practise, all these assumptions are met when
one deals with interaction potential $V(s)=1/r^s$ in dimension $d=3$, 
with $s\in (3,\infty)$.

\subsection{Goals, existing results and difficulties}

We study in this paper the well-posedness of the spatially 
homogeneous Boltzmann equation for singular collision kernel as 
introduced above. In particular we focus on the questions of uniqueness 
and stability with respect to the initial condition which were open, for 
collision kernel with angular cutoff, until the two recent papers \cite{f,dm} 
(except in the special case of Maxwell molecules, see below).

\vip 

In the case of a collision kernel with angular cutoff, that is when 
$\int_0^\pi \beta(d \theta) <+ \infty$, there are some optimal existence and 
uniqueness results: Mischler-Wennberg \cite{mw} in the space of $L^1$ 
non-negative functions with 
finite non-increasing kinetic energy (for counter-examples of spurious 
solutions 
with increasing kinetic energy, see \cite{wenn}  in the hard spheres 
case, and 
\cite{luwenn}  in the case of hard potentials with or without angular 
cutoff), 
Lu-Mouhot \cite{lumou} in the space of non-negative measures with finite 
non-increasing kinetic 
energy.  

\vip

However, the case of collision kernels without cutoff is much more 
difficult. At 
the same time it is crucial from the physical viewpoint since it corresponds 
to the fundamental class of the interactions deriving from inverse power-law 
between particles.  This difficulty is not surprising, since there is 
a difference 
of nature in the collision process between the two cases:
on each compact time interval, each particle collides with infinitely 
(resp. finitely) many others 
in the case without (resp. with) cutoff.

\vip

Until recently, the only uniqueness result obtained 
for non cutoff 
collision kernel was concerning Maxwellian molecules, studied
successively by Tanaka \cite{tanaka}, Horowitz-Karandikar \cite{kh},
Toscani-Villani \cite{tv}: it was proved in \cite{tv} 
that uniqueness holds for the Boltzmann equation 
as soon as
$\Phi$ is constant and (\ref{gccs}) is met, for any 
initial (measure) datum with finite mass and energy, that is
$\int_{\rd} (1+|v|^2)\, f_0(dv)< +\infty$.

\vip

There has been recently two papers in the case where $\beta$ is non 
cutoff and 
$\Phi$ is not constant. The case where $\Phi$ is bounded (together 
with additionnal 
regularity assumptions) was treated in Fournier \cite{f}, for 
essentially any initial (measure) datum such that
$\int_{\rr^d} (1+|v|)f_0(dv)<\infty$. 
More realistic collision kernels have been treated 
by Desvillettes-Mouhot \cite{dm} (including the physical important cases 
of hard and moderately soft potentials without cutoff), for initial 
data in some weighted 
$W^{1,1}$ spaces.

\vip

In the present paper, we extend and improve the method of \cite{f}: 
\begin{itemize} 
\item it can deal with the physical collision kernels corresponding 
to hard and moderately soft 
potentials, as in \cite{dm}: in dimension $d=3$ we obtain well-posedeness
for interaction potentials $1/r^s$ with $s \in (3,\infty)$,
\item the proof is simplified as compared to \cite{f}: it is shorter,
allows measure initial conditions (for technical reasons, we had
to consider only functions in \cite{f}), and 
it does not refer anymore to probabilistic arguments.
\end{itemize} 

Finally let us compare our results with those in \cite{dm},
when applied to the case of an interaction potential $V(s)=1/r^s$
in dimension $d=3$. 
\begin{itemize}
\item Our result is much better in the case of moderately soft potentials
($s \in (3,5)$). Indeed,
we assume only that the initial condition $f_0$ has finite mass, 
energy and entropy
(plus, if $s \in (3,3.48)$, a moment condition 
$\int_{\R^d} |v|^q f_0(v)dv<\infty$ for $q$ large enough).
All these conditions, together with $f_0 \in L^p(\R^d)\cap W^{1,1}(\R^d,
(1+|v|)^2dv)$ (for some $p>1$ depending on the collision rate) 
were assumed in \cite{dm}.  
\item Our result is different in the case of hard potentials 
($s\in (5,\infty)$).
We allow any measure initial $f_0$ condition such that
for some $\e>0$,
$\int_{\R^d} e^{\e |v|^{\gamma}}f_0(dv)<\infty$, where $\gamma=(s-5)/(s-1)$.
In \cite{dm}, the case where $f_0 \in W^{1,1}(\R^d, (1+|v|)^2dv)$
was treated. We thus assume much less regularity, but 
much more {\it localization}.
\end{itemize}

Let us remark that our result is quasi-optimal when $s\in (3.48,5)$,
since the finiteness of entropy and energy is physically
very reasonnable. It might be possible to relax the
entropy condition, but it is not clear: one reasonnably has
to assume a few regularity on $f_0$ to get the uniqueness, since 
the collision rate involves $|v-v_*|^{\gamma}$ with $\gamma<0$,
and we remark that $|v-v_*|^{\gamma}f_0(dv)f_0(dv_*)$ is  
infinite when $f_0$ contains, e.g., Dirac measures.

\vip

Let us emphasize that, as in \cite{dm,f}, we are only able 
to prove well-posedness in the case of a moderate angular singularity 
(assumption ${\bf (A2)}$). 

To our knowledge, there is no uniqueness result under the general assumption
(\ref{gccs}),
except for Maxwellian molecules (see \cite{tv}).

\subsection{Notation} 
Let us denote by $\lip$ the set of globally Lipschitz functions 
$\varphi:\rd \mapsto \rr$, and by $\lipu$ the set of functions 
$\varphi \in \lip$ such that 
$$
\| \varphi \|_{\mbox{{\scriptsize Lip}}(\mathbb{R}^d)} = \sup_{v\ne\tv} 
\frac{|\varphi(v)-\varphi(\tv)|}{|v-\tv|} \leq 1.
$$
Let also $L^p(\rd)$ denote the Lebesgue space of measurable 
functions $f$ such that
$$
\| f \|_{L^p(\rd)} := \left( \int_{\rd} f^p \, dv \right)^{1/p} < +\infty. 
$$

\vip 

Let $\cp(\rd)$ be the set of probability measures on 
$\rd$, and 
\begin{equation*}
\pun = \left\{f \in \cp(\rd), \; m_1(f) <\infty \right\} \quad 
\mbox{ with } \quad  
m_1(f) := \intrd |v| \, f(dv).
\end{equation*}
We denote by $\li$ the set of measurable families $(f_t)_{t\in [0,T]}$ of 
probability measures on $\rd$ such that 
$$
\sup_{[0,T]} \, m_1(f_t) < + \infty,
$$
and by $\lilp$ the set of measurable families $(f_t)_{t\in [0,T]}$ of 
probability measures on $\rd$ such that 
$$
\sup_{[0,T]} \, m_1(f_t) <+ \infty, \quad \int_0 ^T \| f_t \|_{L^p(\rd)} \, dt <+ \infty.
$$

\vip

For $v,v_*\in\rd$, and $\sigma \in \Sp^{d-1}$, we write
\begin{equation*}\label{dfvprime}
v'=v'(v,v_*,\sigma)= \frac{v+v_*}{2} + \frac{|v-v_*|}{2}\sigma,
\end{equation*}
and we write 
$$
\sigma = (\cos \theta, \sin \theta \, \xi) \quad \mbox{ with } \xi \in 
\Sp^{d-2}, \ \theta \in [0,\pi], 
$$ 
in some orthonormal basis of $\R^d$ with first vector $(v-v_*)/|v-v_*|$.

\vip 

Finally we denote $x \land y = \min \{ x,y \}$ and $x_+=\max \{ x,0 \}$, 
and for some 
set $E$ we write $\indiq_E$ the usual indicator function of $E$.

\section{Main results}\label{mainresult} \setcounter{equation}{0}

Let us define the notion of weak (measure) solutions we shall use.

\begin{defin}\label{dfsol}
Let $B$ be a collision kernel which satisfies {\bf (A1-A2)}. 
A family $f=(f_t)_{t\in [0,T]} \in \li$ is a weak solution to (\ref{be}) if
\begin{equation}\label{concon}
\int_0^T dt  \intrd f_t(dv)\intrd f_t(dv_*) \, \Phi(|v-v_*|) \, |v-v_*| 
< + \infty,
\end{equation}
and if for any $\varphi \in \mbox{{\em Lip}}(\rd)$, and any $t\in[0,T]$,
\begin{equation}\label{wbe}
\frac{d}{dt}\intrd \varphi(v)\, f_t(dv) = 
\intrd f_t(dv) \intrd f_t(dv_*) \, A[\varphi](v,v_*),
\end{equation}
where
\begin{equation}\label{afini}
A[\varphi] (v,v_*) = \Phi (|v-v_*|) \, \int_0^\pi \beta(d\theta) 
\int_{\xi \in \Sp^{d-2}} 
\left[\varphi(v')-\varphi(v) \right] \, d\xi.
\end{equation}
\end{defin}

Note that for any $\sigma \in \Sp^{d-1}$,
\begin{equation}\label{majfond} 
|v'-v| = |v-v_*| \, \sqrt{\frac{1-\cost}{2}} \leq  \frac{\theta}{2} \, |v-v_*|,
\end{equation}
so that thanks to assumption {\bf (A2)}, (\ref{concon}) ensures that  
all the terms in (\ref{wbe}) are well-defined.

\vip 

Let us now introduce the distance on $\pun$ we shall use. 
For  $g,\tg \in \pun$, let $\ch(g,\tg)$ be the set of probability 
measures on $\rd\times\rd$ with first marginal $g$ and second marginal $\tg$.
We then set
\begin{eqnarray}\label{dfd1}
d_1(g,\tg)&=&
\inf \left\{\intrdd |v-\tv| \, G(dv,d\tv),\quad G\in \ch(g,\tg)\right\} \ala
&=&\min\left\{\intrdd |v-\tv|\, G(dv,d\tv),\quad G\in \ch(g,\tg)\right\}\ala
&=& \sup \left\{\intrd \varphi(v)\, \big[ g(dv)-\tg(dv) \big], \quad \varphi \in \lipu 
\right\}.
\end{eqnarray}
This distance is the Kantorovitch-Rubinstein distance,
and can be viewed as a particular Wasserstein distance.
We refer to Villani \cite[Section 7]{villani2} for more details on this distance,
and for proofs that the equalities in (\ref{dfd1}) hold.

\vip 

Our main result is the following inequality, which will be applied
in the sequel to hard and soft potentials separately.

\begin{theo}\label{maintheo}
Let $B$ be a collision kernel which satisfies {\bf (A1-A2)}. 
Let us consider two weak solutions  $f,\tf$ 
to (\ref{be}) lying in $\li$, and satisfying
\begin{equation}\label{impo}
\int_0^T dt \int_{\rd \times \rd} \Big[ f_t(dv) \, 
f_t(dv_*)+\tf_t(dv) \, \tf_t(dv_*) \Big]
\, (1+|v|) \, \Phi(|v-v_*|) < +\infty.
\end{equation}
For $s\in [0,T]$, let $R_s\in\ch(f_s,\tf_s)$ be such that 
$$
d_1(f_s,\tf_s)=\int_{\rd\times\rd} |v-\tv| \, R_s(dv,d\tv).
$$ 
Then for all $t \in [0,T]$,
\begin{eqnarray}\label{resultat}
d_1(f_t,\tf_t) \leq d_1(f_0,\tf_0)
+ \kappa_1 \, \frac{|\Sp^{d-2}|}2 \, \intot ds \intrdd R_s(dv,d\tv) 
\intrdd R_s(d\vs,d\tvs) \ala
\times \Big[ 8 \, \big( \Phi(|v-v_*|) \land \Phi(\tv-\tv_*) \big) \, 
|v-\tv| \hskip2cm \ala 
+ \big( \Phi(|v-v_*|) - \Phi(\tv-\tv_*) \big)_+ \, |v-v_*| \hskip1cm\ala
+ \big( \Phi(|\tv-\tv_*|) - \Phi(v-v_*)\big)_+ \, |\tv-\tv_*| \Big].
\end{eqnarray}
\end{theo}

\vip

The meaning of this inequality can be understood 
by means of probabilistic arguments, see \cite{f} for details.
Consider however two infinite particle systems, whose velocity
distributions are $f$ and $\tf$ respectively.
The main ideas are that the first term on the right hand side
expresses an increase of the optimal coupling due to simultaneous
collisions (in both systems), 
whose rate is (optimally) the minimum between the two rates.
Next, the second and third 
terms explain that the optimal coupling also increases
due to a difference between the rates of collision in the two systems. 
Note that these two last terms equal zero in case 
of Maxwellian molecules.

\vip

We now give the application of our inequality to the study
of hard potentials.

\begin{cor}\label{hp} 
Let $B$ be a collision kernel which satisfies {\bf (A1-A2)}, 
and {\bf (A3)}$(\gamma)$ for some $\gamma \in (0,1]$.
\begin{itemize} 
\item[(i)] Let $\e>0$ be fixed. There exists a constant $K_\e>0$ 
such that for any pair
of weak solutions $(f_t)_{t\in[0,T]}$, $(\tf_t)_{t\in[0,T]}$ to 
(\ref{be}), lying in $\li$ and
satisfying
\begin{equation}\label{momex}
C\big(T,f+\tf,\e\big) := \sup_{[0,T]}
\intrd e^{\e |v|^\gamma} \big[f_t+\tf_t \big](dv) < + \infty,
\end{equation}
there holds for all $t \in [0,T]$:
$$
d_1(f_t,\tf_t) \leq d_1(f_0,g_0) + K_\e \, C\big(T,f+\tf,\e\big) 
\, \intot d_1(f_s,\tf_s)
\big(1+ \big|\log d_1(f_s,\tf_s)\big| \big) \, ds.
$$
\item[(ii)] As a consequence for any $f_0 \in \pun$, there exists 
at most one weak solution 
$f \in\li$ to (\ref{be}) starting from $f_0$ and such that $C(T,f,\e)<+\infty$.
\item[(iii)] Let us now give an existence and uniqueness result,
assuming (here only) additionnally {\bf (A4)}$(\gamma)$ and ${\bf (A5)}$. 
Consider $f_0 \in \pun$ such that, for some $\e_0 >0$, $K>0$, we have  
\begin{equation}\label{expini}
\int_\rd e^{\e_0 |v|^\gamma} \, f_0(dv) \leq K < +\infty.
\end{equation}
Then there exists a unique weak solution 
$(f_t)_{t\in[0,\infty)} \in L^1_{\mbox{{\scriptsize {\em loc}}}}([0,\infty), \pun)$ starting from $f_0$.
Furthermore, there exist $\e_1>0$ and $\bar K>0$, depending only on 
$\e_0,K,B$, such that for all $T>0$, $C(T,f,\e_1) \leq \bar K$. 
\item[(iv)] Finally let us give a result on the dependence according to 
the initial datum. Consider a family $(f^n)_{n \geq 1},f^\infty$ 
of weak solutions to (\ref{be}) such that, for some $\e>0$, $T>0$,
we have 
$$
\sup_{n \ge 1} \, C(T,f^\infty+f^n,\e) <+ \infty.
$$
Then 
$$
\lim_{n \to \infty} \, d_1(f^n_0,f^\infty_0)=0
\quad \Longrightarrow \quad 
\lim_{n \to \infty} \sup_{[0,T]} \, d_1(f^n_t,f^\infty_t)=0.
$$
\end{itemize}
\end{cor}

Let us recall that this result applies in particular to hard potentials 
in dimension $d=3$ (that is inverse power-law potentials with $s > 5$). 
In \cite{dm}, under very similar conditions on the collision kernel, 
a well-posedness and stability result was obtained 
in the space $L^\infty([0,T], W^{1,1}(\rd, (1+|v|^2) \, dv))$.
We thus relax the regularity assumption, but we require more moments.

\vip

We finally apply our inequality to the study of soft potentials.

\begin{cor}\label{sp}
Let $B$ be a collision kernel which satisfies {\bf (A1-A2)},
and {\bf (A3)}$(\gamma)$ for some $\gamma \in (-d,0)$. 

\begin{itemize} 
\item[(i)] Let $p\in (d/(d+\gamma),\infty]$ be fixed.
There exists a constant $K_p>0$ such that for any pair
of weak solutions 
$(f_t)_{t\in[0,T]}$, $(\tf_t)_{t\in[0,T]}$ to (\ref{be}) on $[0,T]$, 
lying in $L^\infty([0,T],\pun) \cap L^1([0,T], L^p(\rd))$, there holds
\begin{equation*}
\forall \, t \in [0,T], \qquad d_1(f_t,\tf_t) \leq d_1(f_0,g_0) \, 
e^{K_p \big[C(t,f,p)+C(t,\tf,p)+ t \big]},
\end{equation*}
where 
$$
\forall \, t \in [0,T], \qquad C(t,f,p)=\int_0 ^t \| f_s \|_{L^p(\rd)} \, ds.
$$
Uniqueness and stability thus 
hold in $L^\infty([0,T],\pun) \cap L^1([0,T], L^p(\rd))$.
\item[(ii)] Let $p\in(d/(d+\gamma),\infty]$. For any initial condition 
$f_0 \in \pun \cap L^p(\rd)$, local existence and uniqueness hold, that is 
there exists  
$$
T_*=T_* \big(\|f_0\|_{L^p(\rd)},B \big)>0
$$ 
such that there exists a unique weak solution $(f_t)_{t \in [0,T_*)}$ 
to (\ref{be}) 
which furthermore belongs to
$$
L^\infty_{\mbox{{\scriptsize {\em loc}}}}\big([0,T_*),\pun\cap L^p(\rd)\big).
$$
\item[(iii)] Assume now furthermore that $\gamma \in (-1,0)$, 
{\bf (A4)}$(\gamma)$, and {\bf (A6)}$(\nu)$
for some $\nu \in (-\gamma,1)$. 
Consider an initial datum $f_0\in \pun$ 
with finite energy and entropy, that is
\begin{equation}\label{fee}
\intrd f_0(v) (|v|^2 + |\log f_0(v)| )dv <\infty.
\end{equation}
Assume also that for some $q>\gamma^2/(\nu+\gamma)$, $f_0\in L^1(\R^d,
|v|^q dv)$.
Then there exists a unique weak solution $(f_t)_{t\in [0,\infty)}$
to (\ref{be}), which furthermore belongs to
$$
L^\infty_{\mbox{{\scriptsize {\em loc}}}} \big([0,\infty),\pun  
\cap L^1(\rd,(|v|^q+|v|^2 )\, dv)\big) 
\cap L^1_{\mbox{{\scriptsize {\em loc}}}}  \big([0,\infty), L^p(\rd) \big)
$$
for some (explicit) $p\in (d/(d+\gamma),d/(d-\nu))$.
\end{itemize}
\end{cor}

Let us recall that point (iii) applies, in dimension $d=3$, 
to the case of moderately
soft potentials, that is inverse power-law potentials with 
$s \in (3,5)$. In such a case,
one has $\gamma=(s-5)/(s-1)$
and $\nu=2/(s-1) \in (-\gamma,1)$. We observe that for $s \in (s_0,5)$,
with $s_0=2\sqrt{5}-1\simeq 3.472$, 
the choice $q=2$ is possible,
so that our conditions reduce to the finiteness of entropy and energy.

On the contrary, for $s>3$ close to $3$, $q$ has to be chosen very large,
e.g., for $s=3.01$, we have to take $q\simeq 200$.

\vip

A similar result was obtained in 
\cite[Theorem 1.3]{dm}, assuming that 
$f_0 \in L^p(\rd) \cap  L^1(\rd,|v|^q \, dv)  
\cap W^{1,1}(\rd,(1+|v|^2)\, dv)$, with $p>d/(d+\gamma)$
and $q>\gamma^2/(\nu+\gamma)$.
We thus relax a large part of these conditions.

\vip

The rest of the paper is dedicated to the proof of these results:
we establish Theorem \ref{maintheo} in Section \ref{pr1}.
Applications to hard and soft potentials are
studied in Sections \ref{pr2} and \ref{pr3} respectively.

\section{The general inequality}\label{pr1}\setcounter{equation}{0}

As a preliminary step, we shall parameterize precisely 
the post-collisional velocities.
We follow here 
the approach of \cite{fm}, which was strongly
inspired by Tanaka \cite{tanaka}, and we extend it 
to any dimension $d \ge 2$. 

\vip

The first step is to define a parameterization of the sphere orthogonal  
to some given vector $X \in \R^d$. This parameterization shall 
not be smooth of course. We identify in the sequel $\Sp^0 = \{-1,+1\}$. 

For $X\in\rd\backslash \{0\}$, we set $S_X$ to be the  
symmetry with respect to the hyperplane 
$$
H_X = \left( e_d - \frac{X}{|X|} \right)^\bot
$$
(where $e_d = (0, \dots, 0, 1)$) if $e_d \not= X/|X|$, and $S_X = \mbox{Id}$ 
else. We set 
$$
C_X = \left\{ U \in \rd \ ; \ |U| = |X| \mbox{ and } \langle U, 
X \rangle =0 \right\}.
$$

Then we parameterize $C_X$ by $\Sp^{d-2}$ as follows: we set 
$$
\forall \, \xi=(\xi_1,...,\xi_{d-1}) \in \Sp^{d-2}, \quad \Pi(\xi) = 
(\xi_1, \dots, \xi_{d-1},0) \in \Sp^{d-1} \subset \rd
$$
and 
$$
\Gamma(X,\xi) = |X| \, S_X \big(\Pi(\xi)\big). 
$$

It is easy to check that for a given $X$, the map $\xi \in \Sp^{d-2} 
\mapsto \Gamma(X,\xi)$ is a bijection 
onto $C_X$  and is a unitary parameterization. Therefore, 
for $\xi \in \Sp^{d-2}$, $\theta \in [0,\pi]$, and $X,v,v_*\in \rd$, 
one may write 
$$
v'=v'(v,v_*,\theta,\xi)=v+ \frac{\cost-1}{2} \, 
(v-v_*)+\frac{\sint}{2} \, \Gamma(v-v_*,\xi)
$$
and for all $\varphi\in \mbox{Lip}(\rd)$, recalling (\ref{afini})
$$
A[\varphi](v,v_*)=
\Phi(|v-v_*|) \, \int_0^\pi \beta(d\theta)
\int_{\Sp^{d-2}} d\xi \, \big[ \varphi \big(v'(v,v_*,\theta,\xi)\big)-
\varphi(v) \big].
$$

A problem of this parameterization is its lack of smoothness. To overcome
this difficulty, we shall prove the following fine version of a Lemma
due to Tanaka \cite{tanaka}, whose proof may be found in \cite[Lemma 2.6]{fm} 
in dimension $3$. 

\begin{lem}\label{tanana}
There exists a measurable map
$\xi_0 : \rdd \times \Sp^{d-2} \mapsto \Sp^{d-2}$ 
such that for any $X,Y \in \rd \setminus \{ 0 \}$, 
the map $\xi \mapsto \xi_0(X,Y,\xi)$ is a bijection with jacobian $1$ 
from $\Sp^{d-2}$ into itself (when $d \ge 3$), and 
\begin{equation}\label{controlgamma}
\forall \, \xi \in \Sp^{d-2}, \quad 
\left|\Gamma(X,\xi)- \Gamma\big(Y,\xi_0(X,Y,\xi)\big)\right| \leq 3 \, | X-Y |.
\end{equation}

This implies that for all $v,v_*,\tv,\tv_* \in \rd$, all $\theta\in[0,\pi]$,
all $\xi \in \Sp^{d-2}$, we have 
\begin{eqnarray}\label{delta}
&\big|v'(v,v_*,\theta,\xi)-v'(\tv,\tv_*,\theta,\xi_0(v-v_*,\tv-\tv_*,\xi)\big|
 \ala
& \hspace{5cm} \leq |v-\tv|+ 2 \, \theta \, \big(|v-\tv|+|v_*-\tv_*| \big).
\end{eqnarray}
\end{lem}

\begin{proof}[Proof of Lemma~\ref{tanana}] 
The case $d = 2$ is trivial, therefore we assume $d \ge 3$. 

Let us consider $X,Y \in  \rd \setminus \{ 0 \}$. If $X/|X| = Y/|Y|$, 
it is enough to choose $\xi_0(X,Y,\xi) = \xi$. Indeed in this case 
$S_X=S_Y$ so that
$$
\big| \Gamma(X,\xi) - \Gamma(Y,\xi) \big| = \big| |X| - |Y| \big| \le 
| X- Y |. 
$$
Now assume that $X/|X| \not= Y/|Y|$. Then let us define $R_{X,Y}$ to be the 
axial rotation of $\rd$ transforming $X/|X|$ into $Y/|Y|$ around a line perpendicular 
to the plane determined by $X$ and $Y$. Let us then define $\xi_0$ by the identity 
$$
\Gamma\big(Y, \xi_0(X,Y,\xi) \big) = \frac{|Y|}{|X|} \, R_{X,Y} 
\big( \Gamma(X,\xi) \big) \in C_Y.
$$

For any $X,Y \in \rd \setminus \{0\}$, the application 
$\xi \mapsto \xi_0(X,Y,\xi)$ is the restriction 
to $\Sp^{d-2}$ of the following orthogonal linear transformation on 
$\mathbb{R}^{d-1}$ 
$$
\forall \, Z \in \mathbb{R}^{d-1}, \quad O_{X,Y} (Z) = \Pi^{-1} \circ 
S_Y \circ R_{X,Y} \circ S_X \circ \Pi (Z). 
$$
Therefore it has unit jacobian. Finally let us check the control 
\eqref{controlgamma}: 
\begin{eqnarray*} 
\left| \Gamma(X,\xi) - \Gamma \big( Y, \xi_0(X,Y, \xi ) \big) \right| = 
 \left| \Gamma(X,\xi) - \frac{|Y|}{|X|} \, R_{X,Y} \, \Gamma ( X, \xi ) 
\right| \\ 
 \le \left| \Gamma(X,\xi) \, \left( 1 - \frac{|Y|}{|X|} \right) \right| 
 + \frac{|Y|}{|X|} \,  \left| \Gamma(X,\xi) - R_{X,Y} \, \Gamma ( X, \xi ) 
\big) \right|  \\ 
  \le \left| X - Y \right| 
 + |Y| \,  \left|  \frac{Y}{|Y|} -  \frac{X}{|X|}\right| \le 3 \, |X-Y|.
\end{eqnarray*}
\end{proof} 

Since the transformaction $\xi_0(X,Y,\cdot)$ has unit jacobian, one may 
finally rewrite (\ref{afini}),
for all $\varphi\in \mbox{Lip}(\rd)$, all $X,Y \in \R^d$ (which may 
depend on $v,v_*,\theta$), as
\begin{equation}\label{aphi2}
A[\varphi](v,v_*)=
\Phi(|v-v_*|) \, \int_0^\pi \beta(d\theta)
\int_{\Sp^{d-2}} d\xi \, \big[ \varphi \big(v'\big(v,v_*,\theta,
\xi_0(X,Y,\xi) \big)\big)-\varphi(v) \big].
\end{equation}

\vip

We may finally give the 

\vip

\begin{proof}[Proof of Theorem \ref{maintheo}.] 
We denote 
$$
h^{\varphi}_t:= \int_{\rd} \varphi(v) \, \big(f_t-\tf_t \big)(dv)
$$ 
for $\varphi \in \lipu$, $t\in [0,T]$. 
We also set $h_t=d_1(f_t,\tf_t)$, and we recall that
$$
h_t=\intrdd|v-\tv|\, R_t(dv,d\tv) = \sup_{\varphi \, \in \, 
\mbox{{\scriptsize Lip}}_1(\rd)} h_t^\varphi.
$$

{\it Step 1}.
Let us thus consider $\varphi \in \lipu$.
Using (\ref{wbe}), that $R_t \in \ch(f_t,\tf_t)$ and (\ref{aphi2}),
we immediately obtain, using the map $\xi_0$ built in Lemma \ref{tanana},
\begin{eqnarray}\label{dh1}
\frac{d}{dt}h^\varphi_t &=& \intrdd f_t(dv) \, f_t(dv_*) \, A[\varphi](v,v_*) 
-  \intrdd \tf_t(d\tv) \, f_t(d\tvs) \, A[\varphi](\tv,\tvs) \ala
&=& \intrdd R_t(dv,d\tv) \intrdd R_t(dv_*,d\tvs) \, \Big(A[\varphi](v,v_*)
-A[\varphi](\tv,\tvs) \Big)\ala
&=&  \intrdd R_t(dv,d\tv)
\intrdd R_t(dv_*,d\tvs) \intzp \beta(d\theta)\int_{\Sp^{d-2}} d\xi  \ala
&&\hskip1cm  \Big( \Phi(|v-v_*|) \, 
\big[ \varphi\big( v'(v,v_*,\theta,\xi)\big)-\varphi(v)\big] \\
&&\hskip1cm  -\Phi(|\tv-\tv_*|) \, 
\big[\varphi\big(v'(\tv,\tv_*,\theta,\xi_0 \big(v-v_*,\tv-\tv_*,\xi\big)\big)
-\varphi(\tv)\big] \Big). \nonumber
\end{eqnarray}
We now use the shortened notation 
$$
v'=v'(v,v_*,\theta,\xi) \ \mbox{ and } \ 
\tv'=v'\big(\tv,\tv_*,\theta,\xi_0\big(v-v_*,\tv-\tv_*,\xi \big)\big).
$$
Noting that for all $x,y \in \rr$, $x=x\land y + (x-y)_+$, we easily 
deduce from (\ref{dh1}) that 
\begin{eqnarray*}
\frac{d}{dt}h^\varphi_t &=&\intrdd R_t(dv,d\tv)
\intrdd R_t(dv_*,d\tvs)\intzp \beta(d\theta)\int_{\Sp^{d-2}} d\xi  \ala
&& \hskip1.4cm \Big( \big[ \Phi(|v-v_*|)\land \Phi(|\tv-\tvs|)\big] \times
\big[\varphi(v')-\varphi(\tv') - \varphi(v)+\varphi(\tv)\big] \ala
&&\hskip1.7cm + \big[ \Phi(|v-v_*|) - \Phi(|\tv-\tvs|)\big]_+ \times 
\big[\varphi(v')-\varphi(v)\big] \ala
&&\hskip1.7cm + \big[ \Phi(|\tv-\tv_*|) - \Phi(|v-v_*|)\big]_+ \times
\big[\varphi(\tv)-\varphi(\tv')\big] \Big)\ala
&=:& I^\varphi_1(t)+I^\varphi_2(t)+I^\varphi_3(t),
\end{eqnarray*}
where the last equality stands for a definition. 
Using that $\varphi \in \lipu$, (\ref{majfond}), and {\bf (A2)}, we get
\begin{multline}\label{dh3}
I^\varphi_2(t) +I^\varphi_3(t) \leq  
\kappa_1 \, \frac{|\Sp^{d-2}|}2 \, \intrdd R_t(dv,d\tv)
\intrdd R_t(dv_*,d\tvs)  \ala 
\Big( \big[ \Phi(|v-v_*|) - \Phi(|\tv-\tvs|) \big]_+ \, |v-v_*| \ala 
+ \big[ \Phi(|\tv-\tv_*|) - \Phi(|v-v_*|) \big]_+\, |\tv-\tv_*|\Big).
\end{multline}
Next, using again that $\varphi \in \lipu$, we get that
for all $\e \in (0,\pi)$,
\begin{eqnarray*}
I^\varphi_1(t) & \leq & \intrdd R_t(dv,d\tv)
\intrdd R_t(dv_*,d\tvs)\int_0^\e \beta(d\theta)\int_{\Sp^{d-2}} d\xi  \ala
&& \hskip1cm \big[ \Phi(|v-v_*|) \land \Phi(|\tv-\tvs|) \big] \times
\big[ |v'-v| + |\tv'-\tv| \big] \ala
&&+\intrdd R_t(dv,d\tv)
\intrdd R_t(dv_*,d\tvs)\int_\e^\pi \beta(d\theta) \int_{\Sp^{d-2}} d\xi  \ala
&& \hskip1cm \big[ \Phi(|v-v_*|)\land \Phi(|\tv-\tvs|) \big] \times
\big[|v'-\tv'| - |v-\tv|\big] \ala
&&+\intrdd R_t(dv,d\tv)
\intrdd R_t(dv_*,d\tvs)\int_\e^\pi \beta(d\theta)\int_{\Sp^{d-2}} d\xi \ala
&& \hskip1cm \big[ \Phi(|v-v_*|)\land \Phi(|\tv-\tvs|)\big] \times
\big[|v-\tv| - \big(\varphi(v)-\varphi(\tv)\big) \big] \ala 
&=:& J_1^\varphi(t,\e)+J_2^\varphi(t,\e)+J_3^\varphi(t,\e),
\end{eqnarray*}
where the last equality stands for a definition.
First for $J_2 ^\varphi(t,\e)$, using (\ref{delta}) and {\bf (A2)}, 
we immediately get, by symmetry, that
\begin{eqnarray*}
J_2^\varphi(t) &\leq& 2 \, \left( \int_\e^\pi \theta \, \beta(d\theta)
\right) |\Sp^{d-2}| \, \intrdd R_t(dv,d\tv)
\intrdd R_t(dv_*,d\tvs)\ala
&& \hskip1cm \big[ \Phi(|v-v_*|) \land \Phi(|\tv-\tvs|) \big] \,
\big[ |v-\tv| + |v_*-\tv_*| \big] \ala
&\leq& 4 \, \kappa_1 \,  |\Sp^{d-2}| \, \intrdd R_t(dv,d\tv) 
\intrdd R_t(dv_*,d\tvs) \ala
&& \hskip1cm \big[ \Phi(|v-v_*|) \land \Phi(|\tv-\tvs|)\big] \, |v-\tv|.
\end{eqnarray*}

Next, setting 
$$
\alpha_\e= |\Sp^{d-2}| \, \int_0^\e \theta \, \beta(d\theta), 
$$
it is not hard to obtain, using (\ref{majfond}), 
the fact that $R_t \in \ch(f_t,\tf_t)$ and a symmetry argument, that
\begin{multline}
J^\varphi_1(t,\e)\leq  \frac{\alpha_\e}{2} \, \intrdd R_t(dv,d\tv)
\intrdd R_t(dv_*,d\tvs)\ala 
\times \big[ \Phi(|v-v_*|)\land \Phi(|\tv-\tvs|) \big] \,  
\big( |v|+|\tv|+|v_*|+|\tvs| \big)\ala
\leq \alpha_\e \, \intrd f_t(dv) \intrd f_t(dv_*) \, \Phi(|v-v_*|) \,  |v|
+  \alpha_\e \, \intrd \tf_t(d\tv)\intrd \tf_t(d\tv_*) \, 
\Phi(|\tv-\tv_*|) \, |\tv|\ala
\leq C\big(t,f,\tf \big) \, \alpha_\e, \qquad
\end{multline}
where the constant $C\big(t,f, \tf\big)$ belongs to $L^1([0,T])$ 
due to (\ref{impo}).

\vip 

Finally for $J_3^\varphi(t,\e)$ we notice that the integrand is 
nonnegative (since $\varphi \in \lipu$) and does not depend on 
$\theta,\varphi$.
Hence, denoting 
$$
S_\e:= |\Sp^{d-2}| \, \int_\e^\pi \beta(d\theta) < +\infty, 
$$
we have, for any $A>0$, 
$$
J^\varphi_3(t,\e) \leq K^\varphi_1(t,\e,A)
+ K^\varphi_2(t,\e,A), 
$$
where
\begin{eqnarray*}
K^\varphi_1(t,\e,A) &= & A \, S_\e \, \intrdd R_t(dv,d\tv)
\intrdd R_t(dv_*,d\tvs)\ala
&& \hskip3cm  \big[ |v-\tv|- \big(\varphi(v)-\varphi(\tv)\big) \big]\ala
K^\varphi_2(t,\e,A) &=& S_\e \intrdd R_t(dv,d\tv)
\intrdd R_t(dv_*,d\tvs)\ala
&& \big[ \Phi(|v-v_*|)\land \Phi(|\tv-\tvs|) \big] \, 
\indiq_{\{ \Phi(|v-v_*|)\land \Phi(|\tv-\tvs|) >A\}} \,  |v-\tv|.
\end{eqnarray*}
Using that $R_t \in \ch(f_t,\tf_t)$, and that it achieves the Wasserstein 
distance, we get 
$$
K^\varphi_1(t,\e,A) =  A \, S_\e \, \left[d_1(f_t,\tf_t)-h_t^\varphi \right].
$$
Next, we obtain
\begin{eqnarray*}
K^\varphi_2(t,\e,A)&\leq & S_\e \, \intrd f_t(dv) \intrd f_t (dv_*) \, 
|v| \, \Phi(|v-v_*|) \, \indiq_{\{\Phi(|v-v_*|)>A\}} \ala
&& + S_\e \, \intrd \tf_t(d\tv) \intrd \tf_t (d\tv_*) \, 
|\tv| \, \Phi(|\tv-\tv_*|) \, \indiq_{\{\Phi(|\tv-\tv_*|)>A\}} \ala
&\leq& S_\e \, C_A\big(t,f,\tf\big).
\end{eqnarray*}
Due to (\ref{impo}), we observe that 
$$
\lim_{A\to\infty}\int_0^T C_A\big(t,f,\tf \big) \, dt=0.
$$

\vip

{\it Step 2.} Gathering all the previous estimates, we observe that
for any $\varphi \in \lipu$, $t\in [0,T]$, $\e>0$, $A>0$, we have 
\begin{eqnarray}\label{inedi}
\frac{d}{dt} h_t^\varphi &\leq & H_t + \Gamma_{\e,A}(t)
+ A \, S_\e \, \big[h_t - h^\varphi_t \big],
\end{eqnarray}
where 
$$
\Gamma_{\e,A}(t):=\alpha_\e \, C\big(t,f,\tf\big) + S_\e \, 
C_A\big(t,f,\tf\big), 
$$
and
\begin{multline}
H_t := \kappa_1 \, \frac{|\Sp^{d-2}|}2 \, \intrdd R_t(dv,d\tv) 
\intrdd R_t(dv_*,d\tvs) \Big( 8 \, \big[ \Phi(|v-v_*|) \land 
\Phi(\tv-\tv_*)\big] \, |v-\tv|\ala 
+ \big[ \Phi(|v-v_*|) - \Phi(\tv-\tv_*) \big]_+ \, |v-v_*|
+ \big[ \Phi(|\tv-\tv_*|) - \Phi(v-v_*) \big]_+ \, |\tv-\tv_*|
\Big).
\end{multline}
Recall that $h_t=\sup_{\varphi\in\lipu}h_t^\varphi$, and that our 
aim is to prove that
\begin{eqnarray}\label{cqv}
h_t \leq h_0 + \intot H_s \, ds.
\end{eqnarray}
We immediately deduce from (\ref{inedi}) that
\begin{eqnarray*}
h_t^\varphi \, e^{A \, S_\e \, t} \leq h_0^\varphi + \intot  
e^{A \, S_\e \, s} \, 
\big[ H_s + \Gamma_{\e,A}(s) \big] \, ds 
+ A \, S_\e \,\intot h_s \, e^{A \, S_\e \, s} \, ds.
\end{eqnarray*}
Then we take the supremum over $\varphi\in\lipu$ and we use the 
generalized Gronwall Lemma which states that 
$$
u_t \leq g_t + a\int_0^t u_s \, ds
$$
implies that 
$$
u_t \leq g_0 \, e^{at}+ \int_0^t e^{a\, (t-s)} \, \frac{d g_s}{ds} \,  ds,
$$
which yields
\begin{eqnarray*}
h_t \, e^{A \, S_\e \, t} \leq h_0 \, e^{A \, S_\e \, t} + 
e^{A \, S_\e \, t} \intot
\big[ H_s + \Gamma_{\e,A}(s) \big] \, ds,
\end{eqnarray*}
so that for all $t\in[0,T]$,
\begin{eqnarray*}
h_t \leq h_0 + \intot H_s \, ds + \int_0^T \Gamma_{\e,A}(t) \, dt.
\end{eqnarray*}
This inequality holding for any $\e>0$, $A>0$, we easily conclude that
(\ref{cqv}) holds, since 
$$
\lim_{A\to\infty} \int_0^T \Gamma_{\e,A}(t) \, dt
=\alpha_\e \, \int_0^T C\big(t,f,\tf \big) \, dt
$$ 
with 
$$
\int_0^T C \big(t,f,\tf \big) \, dt < +\infty
$$
and 
$$
\alpha_\e= | \Sp^{d-2} | \, \int_0^\e \theta \, \beta(d\theta) 
\xrightarrow[\e \to 0]{} 0
$$ 
due to {\bf (A2)}.
\end{proof}

\section{Application to hard potentials}
\label{pr2} \setcounter{equation}{0}

\subsection{Propagation of exponential moments} 
We first prove a lemma on the propagation (and appearance)
of exponential moment, which is a variant 
of results first obtained in \cite{Bo97,BoGaPa} (and also developed 
in \cite{gm,MMR1}).

\begin{lem}\label{lem:mts}
Let $B$ be a collision kernel satisfying  
assumptions {\bf (A1-A2-A5)} and {\bf (A3)-(A4)}$(\gamma)$ for some
$\gamma \in (0,1]$.  Let $f_0 \in \pun$.
\begin{itemize}
\item[(i)] Assume that for some $\e_0>0$, some $s\in (0,2)$,
$$
\int_{\rd} e^{\e_0 |v|^s} \, f_0(dv) \le C_{\e_0,s} < +\infty.
$$ 
Then there exists $\e_1>0$ and a constant $C>0$, depending only 
on $s$, $\e_0$, $C_{\e_0,s}$, 
such that for any $T>0$, any weak solution 
$(f_t)_{t\in [0,T]}$ to~\eqref{be} satisfies    
$$
\sup_{[0,T]} \int_{\rd} e^{\e_1 |v|^s} \, f_t(dv) \leq C < +\infty.
$$
\item[(ii)] Assume now only that $e_0=\int_\rd |v|^2 f_0(dv)<\infty$.
For any $s\in (0,\gamma/2)$, any $\tau>0$, there exists $\e>0$ and 
$C>0$, depending only on $s$, $\tau$, and an upperbound of $e_0$
such that for any  $T>0$, any weak solution $(f_t)_{t\in [0,T]}$ to~\eqref{be} satisfies    
$$
\sup_{t\in [\tau,T]} \int_{\rd} e^{\e |v|^s} \, f_t(dv) \leq C < +\infty.
$$
\end{itemize}
\end{lem}

\begin{proof}[Proof of Lemma~\ref{lem:mts}]
We first recall that for any $t\in[0,T]$,
\begin{equation}\label{momencons}
\intrd |v|^2 f_t(dv)=\intrd |v|^2 f_0(dv)=:e_0,
\end{equation}
and we observe that for all $v\in \R^d$, all $t\geq 0$, since 
$\gamma \in (0,1]$ and since $f_t \in \pun$,
\begin{equation}\label{mino}
\intrd |v-v_*|^\gamma f_t(dv_*) \geq |v|^\gamma -
\intrd |v_*|^\gamma f_t(dv_*) \geq   |v|^\gamma - e_0^{\gamma/2}.
\end{equation}
Let us fix $0<s<2$. We define for any $p \in \mathbb{R}_+$
$$ 
m_p(t) := \int_{\rd} |v|^{sp} \, f_t(dv).
$$

{\it Step 1.} The evolution equation (\ref{wbe}) yields 
\begin{equation} \label{mp1}
\frac{d m_p}{dt} =  \int_{\rdd} \Phi(|v-v_*|) \,  
K_p (v,v_*) \, f_t(dv) \, f_t(dv_*),
\end{equation}
where, using {\bf (A5)} and a symmetry argument,
\begin{equation*} 
K_p(v,v_*) := \frac12 \, \int_0^\pi 
\int_{\Sp^{d-2}} \big(|v'|^{sp} + |v'_*|^{sp} - |v|^{sp} - 
|v_*|^{sp}\big) \, b(\cos \theta) \, d\theta\, d\xi.
\end{equation*}
Let us split $b = b^c _\eta + b^r _\eta$ for some $\eta \in (0,\pi)$ 
with 
$$
b^c _\eta (\cos \theta) = b(\cos \theta) \, \indiq_{\theta \ge \eta} 
+ \big[ b(\cos \eta) + b'(\cos \eta) \, (\cos \theta - \cos \eta) \big] \, 
\indiq_{0 \le \theta \le \eta} 
$$
for $\theta \in (0,\pi]$. Due to {\bf (A5)}, we know that $b^c _\eta \leq b$,
so that $b^r _\eta\geq 0$. We can split correspondingly 
$K_p = K_p ^{c,\eta} + K_p ^{r,\eta}$. We also easily check that for each 
$\eta\in (0,\pi)$,  $b^c _\eta$ is convex, non-decreasing, and bounded 
on $[-1,1)$.
We are thus in a position to apply \cite[Corollary 1]{BoGaPa}, which
yields that for $p>2/s$,
\begin{equation} \label{mp3}
K_p ^{c, \eta} (v,v_*) \le \alpha_p(\eta) \, 
\big(|v|^2 + |v_*|^2 \big)^{sp/2} 
- K(\eta) \, \big( |v|^{sp} + |v_*|^{sp} \big)
\end{equation}
where $(\alpha_p(\eta))_p$ is strictly decreasing  and satifies
\begin{equation} \label{mp4}
\forall \, p>2/s, \quad 0< \alpha_p <  \frac{C(\eta)}{sp+1}
\end{equation}
for some constant $C(\eta)$ depending on an upper bound of $b^c_\eta$, 
and some constant $K(\eta)$ depending on a lower bound of the mass of 
$\beta^c _\eta$. 
Therefore $K$ can be made uniform according to $\eta$ as $\eta \to 0$.   

For the other part of the collision kernel we use for instance 
\cite[Lemma~2.1]{dm} 
and assumption {\bf (A2)} to deduce that (as soon as $sp \ge 2$)
\begin{equation} \label{mp5}
K_p ^{r, \eta} (v,v_*) \le \delta(\eta) \, \big(  |v|^{sp} + |v_*|^{sp} \big)
\end{equation}
with
$$
\delta(\eta)\leq \mbox{cst} \, \int_0^\pi \theta \, b^r_\eta(\cos\theta) \, d\theta \to 0
$$
as $\eta \to 0$, due to {\bf (A2)}.
 
Combining (\ref{mp3},\ref{mp4},\ref{mp5}) and fixing carefully 
$\eta$ we thus find for all $p>2/s$
\begin{equation*}
K_p (v,v_*) \le \bar \alpha_p \, \big(|v|^2 + |v_*|^2 \big)^{sp/2} 
- \bar K \, \big( |v|^{sp} + |v_*|^{sp} \big)
\end{equation*}
for some constant $\bar K >0$ and where $(\bar \alpha_p)_p$ is 
strictly decreasing and satisfies, for some constant $\bar C>0$,
\begin{equation*}
\forall \, p>2/s, \quad 0< \bar \alpha_p <  \frac{\bar C}{sp+1}.
\end{equation*}
We of course deduce that for $p$ large enough, say 
$p\geq p_0>2/s$,
\begin{equation} \label{mp6}
K_p (v,v_*) \le \bar \alpha_p \, \left[\big(|v|^2 + |v_*|^2 \big)^{sp/2} 
-|v|^{sp} - |v_*|^{sp} \right]  - \bar K \, \big( |v|^{sp} + |v_*|^{sp} \big),
\end{equation}
changing if necessary the value of $\bar K>0$.

We now insert (\ref{mp6}) in (\ref{mp1}). Using {\bf (A4)}$(\gamma)$ 
and (\ref{mino}), we get, for $p\geq p_0$,
\begin{equation}\label{hereweare}
\frac{d m_p}{dt}
\le \bar \alpha_p \, Q_p - K' \, (m_{p+\gamma/s}-e_0^{\gamma/2}m_p)
\end{equation}
for some new constant $K'>0$ and with
\begin{eqnarray*}
Q_p := \int_{\rdd} \left[ \big(|v|^2 + |v_*|^2 \big)^{sp/2} 
- |v|^{sp} - |v_*|^{sp} \right] \, \Phi(|v-v_*|) \, f_t(dv) \, f_t(dv_*). 
\end{eqnarray*}

{\it Step 2.} Using {\bf (A3)}$(\gamma)$ and 
following line by line the proof of \cite[Lemma 4.7]{gm} 
from \cite[eq. (4.13)]{gm} which is the same as (\ref{hereweare}) here 
to \cite[eq. (4.19)]{gm} (this proof is itself 
essentially based on \cite{BoGaPa}), we obtain the following
conclusion. Set $k_p=[sp/4+1/2]$ (here $[\cdot]$ stands for the integer part).
Set also, with the usual Gamma function,
\[ z_p:=\frac{m_p}{\Gamma(p+1/2)} \quad \hbox{ and } 
\quad 
Z_p := \max_{k=1,..,k_p} \big\{ z_{(2k+\gamma)/s} \, z_{p-2k/s}, \,  
z_{2k/s} \, 
z_{p-2k/s+\gamma/s} \big\} . \]
Then for some constants $A'>0$, $A''>0$, $A''' >0$, for all $p\geq p_0$,
\begin{equation}\label{mp10}
\frac{d z_p}{dt} \le A' \, p^{\gamma/s -1/2} \, Z_p 
- A'' \, p^{\gamma/s} \, z_p^{1+\frac{\gamma}{sp}} +A''' z_p. 
\end{equation}

{\it Step 3.} Next, point (i) can
be checked following the ideas of 
\cite[Proposition~3.2]{MMR1} (for $\gamma =1$) and using that classically,
for any $p\geq 0$, $\sup_{t\in [0,\infty)} m_p \leq C_p$, for
some constant depending only on $B$
and $m_p(0)$, see e.g., \cite[Theorem 1-(ii)]{villani} 
and \cite[Lemma 2.1]{dm}.

\vip

{\it Step 4.} Finally, point (ii) can be proved following line by line
the proof of  \cite[Lemma 4.7]{gm}.

\end{proof}

\subsection{Proof of Corollary \ref{hp}}

We first recall the following variant of a 
classical lemma used by Yudovitch \cite{yudo} 
in his Cauchy theorem for bidimensional 
incompressible non-viscious flow. See \cite[Lemme 5.2.1, p. 89]{che}
for a proof.

\begin{lem}\label{yudovitch}
Consider a nonnegative bounded function $\rho$ on $[0,T]$, 
a real number $a \in [0,\infty)$, and a strictly positive, 
continuous and non-decreasing function $\mu=\mu(x)$ on $(0,\infty)$. 
Assume furthermore that 
$$
\int_0^1 \frac{dx}{\mu(x)} = + \infty,
$$ 
and that for all $t\in [0,T]$, 
$$
\rho(t) \leq a + \int_0^t \mu(\rho(s)) \, ds.
$$

Then
\begin{itemize} 
\item[(i)] if $a=0$, then $\rho(t)= 0$ for all $t\in [0,T]$; 
\vspace{0.2cm} 
\item[(ii)] if $a>0$, then 
$$
\forall \, t \in [0,T], \quad m(a)- m(\rho(t)) \leq t
$$ 
where 
$$
m(x)=\int_x^1 \frac{dy}{\mu(y)}.
$$
\end{itemize}
\end{lem}

We may now give the

\vip

\begin{proof}[Proof of Corollary \ref{hp}.]
We thus consider $\gamma \in (0,1]$, and we assume 
{\bf (A1)-(A2)-(A3)}$(\gamma)$.
We also consider some $\e>0$ fixed. 

\vip

{\it Step 1.} Let us first prove point (i). 
Let us consider two weak solutions $(f_t)_{t\in[0,T]}$, 
$(\tf_t)_{t\in[0,T]}$ to (\ref{be}), 
lying in $\li$ and satisfying (\ref{momex}). 
We are in position to apply Theorem \ref{mainresult},
since {\bf (A3)}$(\gamma)$ and (\ref{momex}) clearly guarantee
that (\ref{impo}) holds.
We thus know that (\ref{resultat}) holds.
Using {\bf (A3)}$(\gamma)$, simple computations show that 
\begin{equation*}
\big( \Phi(|v-v_*|)\land \Phi(|\tv-\tvs|)\big) \, |v-\tv|
\leq C \big[ |v|^\gamma+|v_*|^\gamma \big] \, |v-\tv|,
\end{equation*}
while
\begin{eqnarray*}
&&\big[ \Phi(|v-v_*|)-\Phi(|\tv-\tvs|)\big]_+|v-v_*| \ala
&&\leq C\left| \left(|v-v_*|^\gamma-|\tv-\tvs|^\gamma \right)\right| \,  
\big[|v-v_*|\land |\tv-\tvs|+\left|\left( |v-v_*|-|\tv-\tvs|\right)\right| 
\big]\ala
&& \leq C \gamma \, (|v-v_*|\land |\tv-\tv_*|)^{\gamma-1} \, 
\big| \big( |v-v_*|-|\tv-\tvs| \big) \big| \, (|v-v_*|\land |\tv-\tvs|)   \ala
&& \qquad \qquad \qquad+ C \big[ |v-v_*|^\gamma +|\tv-\tv_*|^\gamma \big] 
\, \big[|v-\tv|+|v_*-\tvs|\big] \ala
&& \leq C(1+\gamma) \big[ |v|^\gamma + |v_*|^\gamma + |\tv|^\gamma + 
|\tv_*|^\gamma \big] \, \big[|v-\tv|+|v_*-\tvs|\big].
\end{eqnarray*}
We hence obtain by inserting these inequalities
in (\ref{resultat}) and using symmetry properties, that 
for some constant $D>0$,
\begin{eqnarray*}
d_1(f_t,\tf_t) \leq d_1(f_0,\tf_0) + D
\intot ds \intrdd R_s(dv,d\tv)
\intrdd R_s(dv_*,d\tvs) \ala 
\mbox{ } \qquad \qquad \times \big[|v|^\gamma+|v_*|^\gamma+|\tv|^\gamma
+|\tv_*|^\gamma \big] \, |v-\tv|.
\end{eqnarray*}
Recall now that $R_s \in \ch(f_s,\tf_s)$ achieves the
Wasserstein distance. It is thus clear (recall that $C(T,f+\tf,\e)$
was defined in (\ref{momex})) that 
$$
\sup_{[0,T]} \int_{\rdd}R_s(dv_*,d\tvs) \, 
\big[|v_*|^\gamma+|\tv_*|^\gamma \big] \leq A_\e \, C(T,f+\tf,\e)
$$ 
for some constant $A_\e$. 
We thus get
\begin{eqnarray*}
d_1(f_t,\tf_t) \leq d_1(f_0,\tf_0) + D \, A_\e \, C\big(T,f+\tf,\e\big)
\intot d_1(f_s,\tf_s) \, ds \ala
+ D \, \intot ds \intrdd R_s(dv,d\tv) \, 
\big[ |v|^\gamma+|\tv|^\gamma \big] \, |v-\tv|.
\end{eqnarray*}

Next, for any $s\in [0,T]$ and $a>0$, we have
\begin{eqnarray*}
&&\intrdd R_s(dv,d\tv) \, \big[ |v|^\gamma + |\tv|^\gamma \big] \, |v-\tv|
\leq 2\, a^\gamma \, d_1(f_s,\tf_s)\ala
&&\hskip1cm+\intrdd R_s(dv,d\tv) 
\big[ |v|^\gamma +|\tv|^\gamma \big] \, \big[ |v|+|\tv| \big] \,  
\big( \indiq_{\{|v|>a\}} +\indiq_{\{|\tv|>a\}} \big)\ala
&&\leq 2 \, a^\gamma \, d_1(f_s,\tf_s) + L_\e \,  \intrdd R_s(dv,d\tv) \, 
\big[ e^{-\e a^\gamma/2} \, 
e^{\e |v|^\gamma} + e^{-\e a^\gamma/2} \, e^{\e |\tv|^\gamma} \big] \ala
&& \leq 2 \, a^\gamma \, d_1(f_s,\tf_s) 
+ L_\e \, e^{-\e a^\gamma/2} \, C\big(T,f+\tf,\e\big),
\end{eqnarray*}
for some constant $L_\e$ such that 
$$
\big[ |v|^\gamma + |\tv|^\gamma \big] \, \big[ |v|+|\tv|\big] \, 
\big[ e^{\e |v|^\gamma/2} + e^{\e |\tv|^\gamma/2} \big]
\leq L_\e \, \big[e^{\e |v|^\gamma}+ e^{\e |\tv|^\gamma} \big].
$$

Choosing $a$ such that 
$$
a^\gamma= \big| 2 \, \log d_1(f_s,\tf_s) /\e \big|,
$$ 
we finally get
\begin{eqnarray*}
\intrdd R_s(dv,d\tv) \, \big[ |v|^\gamma + |\tv|^\gamma \big] \, |v-\tv|
\leq \frac{4}{\e} \, d_1(f_s,\tf_s) \, \big| \log d_1(f_s,\tf_s)\big| \ala
+ L_\e \, C\big(T,f+\tf,\e\big) \, d_1(f_s,\tf_s).
\end{eqnarray*}

We finally obtain, setting 
$K_\e=D \,(A_\e+L_\e+4/\e)$, 
that 
\begin{equation*}
d_1(f_t,\tf_t) \leq d_1(f_0,\tf_0) + K_\e \, C\big(T,f+\tf,\e\big) 
\intot d_1(f_s,\tf_s) \, \big(1+ \big| \log d_1(f_s,\tf_s) \big|\big) \, ds.
\end{equation*}

{\it Step 2.} Points (ii) and (iv) are immediate consequences of point
(i) and Lemma \ref{yudovitch} applied 
with $\mu(x)=x \, \big( 1+|\log x| \big)$.

\vip

{\it Step 3.} Finally, we check point (iii).
We thus assume 
{\bf (A1)-(A2)-(A3)}$(\gamma)${-\bf (A4)}$(\gamma)${-\bf (A5)}, and
consider an initial condition $f_0\in \pun$ satisfying (\ref{expini}) for
some $\e_0>0$. Then we know from Lemma \ref{lem:mts}-(i) that any weak solution
starting from $f_0$ satisfies (\ref{momex}) for some $\e_1>0$.
We thus deduce the uniqueness part from point (ii).

Next, we approximate $f_0$ by a sequence of initial conditions $f_0^n$ 
with finite entropy satisfying (\ref{expini}) uniformly (in $n$), 
and such that $d_1(f_0, f_0^n)$
tends to $0$. Then, using for example the existence
result of Villani \cite[Theorem 1]{villexi}, we know that for each
$n$, there exists a weak solution $(f^n_t)_{t\geq 0}$ 
to (\ref{be}) starting from $f_0^n$.
Due to Lemma \ref{lem:mts}-(i), we deduce that there exists $\e_1>0$ such
that for all $T>0$, $\sup_n C(T,f^n_t,\e_1)<\infty$.
It is then not hard to deduce from point (i) and Lemma \ref{yudovitch} 
that there exists $(f_t)_{t\geq 0}$ such that for all $T>0$,
$C(T,f,\e_1)<\infty$ and  $\lim_n \sup_{[0,T]} d_1(f^n_t,f_t)=0$.
An easy consequence is that $(f_t)_{t\geq 0}$ is a weak solution to (\ref{be})
starting from $f_0$.
\end{proof}

\section{Application to soft potentials}
\label{pr3} \setcounter{equation}{0}

The application to soft potentials is easier, since we shall apply
the standard Gronwall Lemma instead of that of Yudovitch.

\vip

\begin{proof}[Proof of Corollary \ref{sp}.]
We consider $\gamma \in (-d,0)$, and assume that 
{\bf (A1)-(A2)-(A3)}$(\gamma)$.

\vip

We observe at once that for $\alpha \in (-d,0)$, and for
$q\in (d/(d+\alpha),\infty]$, there exists a constant $C_{\alpha,q}$ such that
for any $g \in \pun \cap L^q(\rd)$, any $v\in \R^d$,
\begin{eqnarray}\label{suplq}
\intrd g(v_*) \, |v-v_*|^{\alpha} \, dv_* & \leq & 
\int_{|v_*-v|<1} g(v_*) \, |v-v_*|^{\alpha} \, dv_* + \int_{|v_*-v|\geq 1} 
g(v_*) \, dv_*\ala
&\leq& C_{\alpha,q} \, \| g \|_{L^q(\rd)} +1 .
\end{eqnarray}

{\it Step 1.} We first prove point (i). Let thus $p\in (d/(d+\gamma),\infty]$. 
We consider two solutions $f,\tf$ as in the statement.
In order to apply Theorem \ref{maintheo}, we have to check that 
(\ref{impo}) holds. But using (\ref{suplq}), since $p>d/(d+\gamma)$,
we get for $t \in [0,T]$
\begin{eqnarray*}
&&\intrd f_t(dv) \intrd f_t(dv_*) \, (1+|v|) \, |v-v_*|^\gamma \ala
&&\leq 
\left[\intrd f_t(dv)\, (1+|v|)\right]  \sup_{v \in \rd} 
\intrd f_t(dv_*) \, |v-v_*|^\gamma \ala
&&\leq 
\left[\intrd f_t(dv)\, (1+|v|)\right] \, 
\left( C_{\gamma,p}\| f_t \|_{L^p(\rd)}+1\right).
\end{eqnarray*}
The same estimate holds for $\tf$, and therefore we conclude that the 
estimate (\ref{impo}) holds using that $f$ and $\tf$ belong to $\lilp$. 

\vip

Hence we deduce that (\ref{resultat}) holds.
Simple computations using {\bf (A3)}$(\gamma)$ show that 
\begin{equation*}
\Big( \Phi(|v-v_*|) \land \Phi(|\tv-\tvs|)\Big) \, |v-\tv|
\leq C |v-v_*|^\gamma \, |v-\tv|,
\end{equation*}
while 
\begin{eqnarray*}
&& \Big( \Phi(|v-v_*|) - \Phi(|\tv-\tvs|)\Big)_+ \, |v-v_*| \ala
&& \leq C \Big|\Big(|v-v_*|^\gamma - |\tv-\tvs|^\gamma\Big)\Big| \, 
\Big( |v-v_*|\land |\tv-\tvs| +  
\Big|\Big( |v-v_*| - |\tv-\tvs|\Big)\Big|\Big)\ala
&& \leq C |\gamma| \, \left( |v-v_*|\land |\tv-\tvs|\right)^{\gamma-1} \,
\Big| \left( |v-v_*|-|\tv-\tvs|\right)\Big| \, \left( |v-v_*|\land 
|\tv-\tvs|\right) \ala
&&+ C \Big( |v-v_*|^\gamma \lor |\tv-\tvs|^\gamma\Big) \, 
\Big|\left( |v-v_*|-|\tv-\tvs|\right)\Big|\ala
&& \leq C(1+|\gamma|) \, \left( |v-v_*|^\gamma + |\tv-\tvs|^\gamma \right) \, 
\left(|v-\tv|+|v_*-\tvs| \right).
\end{eqnarray*}

Inserting these inequalities
in (\ref{resultat}) and using a symmetry argument, 
we obtain that for some constant $D>0$,
\begin{eqnarray*}
d_1(f_t,\tf_t) \leq d_1(f_0,\tf_0) +
D \, \intot ds \intrdd R_s(dv,d\tv)\intrdd R_s(dv_*,d\tvs) \ala
\Big[|v-v_*|^\gamma+|\tv-\tv_*|^\gamma \Big] \, |v-\tv|.
\end{eqnarray*}
Recall now that $R_s \in \ch(f_s,\tf_s)$ and achieves the
Wasserstein distance.
Hence,
\begin{eqnarray*}
&&\sup_{v,\tv}  \intrdd R_s(dv_*,d\tvs)
\Big[|v-v_*|^\gamma+|\tv-\tv_*|^\gamma \Big] \ala
&&\leq \sup_{v} \intrd f_t(dv_*) \, |v-v_*|^\gamma +
\sup_{\tv} \intrd \tf_t(d\tv_*) \, |\tv-\tv_*|^\gamma \ala
&&\leq C_{\gamma,p} \, \| f_t \|_{L^p(\rd)} + C_{\gamma,p} 
\, \| \tf_t \|_{L^p(\rd)} 
+2
\ala
&&\leq C_{\gamma,p} \, \big[ \| f_t \|_{L^p(\rd)} +\| \tf_t \|_{L^p(\rd)} \big] +2,
\end{eqnarray*}
where we used (\ref{suplq}). Since
$$
\intrdd R_s(dv,d\tv) \, |v-\tv|=d_1(f_s,\tf_s)
$$
we obtain finally, choosing $K_p:=D \, ( C_{\gamma,p}+2)$,
\begin{eqnarray*}
&d_1(f_t,\tf_t) 
\leq d_1(f_0,\tf_0) + K_p \, \big[\| f_t \|_{L^p(\rd)} + \| \tf_t \|_{L^p(\rd)} +1\big] \, 
\intot d_1(f_s,\tf_s)\, ds.
\end{eqnarray*}
The Gronwall Lemma then allows us to conclude the proof.

\vip

{\it Step 2.} We now check point (ii). 
We only have to prove the existence of solutions, since uniqueness
follows from point (i).
Using some results of Villani \cite[Theorems 1 and 3]{villexi}, 
we know that for $\gamma \in (-d,0)$,
for any $f_0\in \pun$ such that 
$$
\int f_0(v)\, \big(|v|^2+|\log f_0(v)|\big)\, dv < + \infty,
$$ 
there exists a weak solution $f\in L^\infty([0,\infty), (1+|v|^2)\, dv)$ to 
(\ref{be}) starting from $f_0$.
Then the existence result of point (ii) follows immediately from point (i) 
together with the
following {\it a priori} estimates,
which guarantee that if $f_0 \in \pun \cap L^p(\R^d)$, then 
this bound propagates locally (in time):
first there exists $C=C(B)$ such that (see
\cite[Proposition 3.2]{dm} and its proof) for any $p\in (d/(d+\gamma),\infty]$,
any weak solution to (\ref{be}) satisfies
\begin{equation*}
\frac{d}{dt} \| f_t \|_{L^p} \leq C \, \big(1+ \| f_t \|_{L^p(\rd)}^2\big),
\end{equation*} 
so that for $0\leq 
t< T_*:=\frac{1}{C}(\pi/2 - \arctan ||f_0||_{L^p})$, we have 
\begin{equation}\label{tat}
\| f_t \|_{L^p} \leq \tan \left( \arctan \| f_0 \|_{L^p} + C \, t \right).
\end{equation}

\vip

Next, we easily check, using (\ref{wbe}), (\ref{majfond}) and {\bf (A2)} 
that
\begin{eqnarray*}
\frac{d}{dt} \intrd |v| \, f_t(dv) &\leq& \kappa_1 \, \frac{|\Sp^{d-2}|}{2}\, 
\intrd f_t(dv) \intrd f_t(dv_*) \, |v-v_*|^{1+\gamma}.
\end{eqnarray*}
If $1+\gamma \geq 0$, we immediately conclude, since 
$|v-v_*|^{1+\gamma} \leq 1+|v|+|v_*|$, that
\begin{equation*}
\frac{d}{dt} \int_\rd |v| \, f_t(dv) \leq \kappa_1 \, \frac{|\Sp^{d-2}|}{2}\, 
\, \left( 1+2\int_\rd |v| \, f_t(dv) \right),
\end{equation*}
so that for $t\geq 0$, we have 
$$
\int_\rd |v| f_t(dv) \leq e^{\kappa_1|\Sp^{d-2}|t} \, 
\left( \int_\rd |v| f_0(dv) + 1\right).
$$

If $1+\gamma \leq 0$, we use (\ref{suplq}) (with $\alpha=1+\gamma$ and $q=p$,
which is valid since $p>d/(d+\gamma)>d/(d+\alpha)$), and we deduce
that
\begin{equation*}
\frac{d}{dt} \int_\rd |v| \, f_t(dv) \leq \kappa_1 \, \frac{|\Sp^{d-2}|}{2}
\, \left( C_{1+\gamma,p} \, \| f_t \|_{L^p(\rd)} +1 \right)=: A \, \| f_t \|_{L^p(\rd)}+ A',
\end{equation*} 
so that for $0\leq t < T_*$, we have, recalling (\ref{tat}), 
$$
\int_\rd |v| \, f_t(dv) \leq \int_\rd |v| \, f_0(dv) + 
A \, \int_0^t \tan \left( \arctan \| f_0 \|_{L^p} +C\,s \right) \, ds 
+A' \, t.
$$

\vip

{\it Step 3.} 
We now assume additionnally that $\gamma \in (-1,0)$, 
({\bf A4})($\gamma$), and ({\bf A6})($\nu$) for some
$\nu \in (-\gamma,1)$. We consider an initial datum $f_0$
with finite energy and entropy (\ref{fee}), and such that for some
$q>q_0=\gamma^2/(\nu+\gamma)$, $f_0 \in L^1(\R^d, |v|^q dv)$.
Applying the result of Villani \cite[Theorem 1]{villexi}, we know that
there exists a weak solution $(f_t)_{t\in [0,\infty)}$ to (\ref{be}).

To conclude the proof, it suffices to apply point (i), and to 
check that for any weak solution $(f_t)_{t\in [0,\infty)}$ to (\ref{be})
starting from $f_0$,

\vip

(a) $f \in L^\infty_{loc} ([0,\infty), L^1(\R^d,(|v|^2+|v|^q) dv))$, 
\smallskip

(b) there exists $p>p_0:=d/(d+\gamma)$ such that  
$f \in L^1_{loc} ([0,\infty), L^p(\R^d))$.

\vip

Point (a) follows from a straightforward application of (\ref{wbe}),
using ({\bf A1})-({\bf A2})-({\bf A3})($\gamma$) 
and that $\gamma \in (-1,0)$, and concluding with the Gronwall Lemma.

\vip

To check point (b), we follow the line of 
\cite[Proposition 3.3]{dm} 
(see (3.2), (3.3) and (3.4) in \cite{dm}),
which was relying on exploiting the entropy production and its 
regularization property obtained by 
Alexandre-Desvillettes-Villani-Wennberg \cite{advw}. 

Exactly as in \cite[(3.2)]{dm}, we get that for any $\alpha>0$,
\begin{equation}\label{lpe}
\int_0 ^T \| (1+|v|)^{\gamma-\alpha}f_t(v) \|_{L^{d/(d-\nu)}(\rd)} \, dt \le C \, (1+T)
\end{equation}
for any $\alpha >0$ any $T>0$ 
and some constant $C>0$ (depending on $\alpha$). 
Using point (a), we also now that for all $T>0$, 
\begin{equation}\label{lunq}
C_T := \sup_{[0,T]} \| (1+|v|)^q f_t(v)\|_{L^1(\R^d)}<\infty.
\end{equation}

By interpolation between estimates (\ref{lpe}) and (\ref{lunq}), we see that
for all $T>0$, for another constant $C_T$ depending on $T$,
\begin{equation}\label{lunqq}
\int_0 ^T \| f \|_{L^p (\rd)} \, dt \le C_T
\end{equation}  
for any $1 < p < d/(d-\nu)$ as soon as, for instance 
\begin{equation}\label{cqqf}
q > (\alpha - \gamma) \frac{p-1}{1-p(d-\nu)/d}.
\end{equation}
Since $p_0=d/(d+\gamma)<d/(d-\nu)$ (because $\gamma > -\nu$)
and since by assumption,
$$
q>q_0 = (- \gamma) \frac{p_0-1}{1-p_0(d-\nu)/d} = \frac{\gamma^2}{\gamma + \nu},
$$
we clearly have (\ref{cqqf}) when choosing with $\alpha>0$ small enough and
$p>p_0$ close enough. This concludes the proof of point (b).
\end{proof}

\def\refname{References}

\end{document}